\documentclass{article}
\usepackage[T2A]{fontenc}
\usepackage[cp1251]{inputenc}
\usepackage[english,russian]{babel}
\usepackage[tbtags]{amsmath}
\usepackage{amsfonts,amssymb,mathrsfs,amscd}
\numberwithin{equation}{section}
\newtheorem{remark}{Remark}[section]
\newtheorem{lemma}{Lemma}[section]
\newtheorem{theorem}{Theorem}[section]
\newtheorem{definition}{Definition}[section]

\newcommand{\res}{\mathop{\rm res}}

\renewcommand{\Im}{\mathop{\rm Im}}
\renewcommand{\Re}{\mathop{\rm Re}}

\newcommand{\col}{\mathop{\rm col}}
\renewcommand{\d}{\mathop{\rm def}}

\begin{document}
\begin{Large}
\thispagestyle{empty}
\begin{center}
{\bf Direct spectral problem for a Dirac operator with non-local potential\\
\vspace{5mm}
V. A. Zolotarev}\\
\end{center}
\vspace{5mm}

{\small {\bf Abstract.} For a Dirac operator with non-local potential on a finite segment, a method of reconstruction of non-local potential from the spectral data is developed. Description of spectral data for such class of operators is given.}
\vspace{5mm}


{\it Key words}: Dirac operator, non-local potential, inverse problem.
\vspace{5mm}

Unlike the Sturm -- Liouville operator, study of a one-dimensional Dirac operator has its peculiar properties based upon the fact that coefficients of the operator are matrix-valued \cite{1, 2}. This adds some specific features to the solution of inverse problems.

In this work, the selfadjoint Dirac operator on a finite segment with non-local potential
$$(Dy)(x)=\left(
\begin{array}{ccc}
0&1\\
-1&0
\end{array}\right)y'(x)+\mu\int\limits_0^\pi(y_1(t)\overline{v_1(t)}+y_2(t)\overline{v_2(t)})dtv(x)$$
is studied where $\mu\in\mathbb{R}$, $v(x)=\col[v_1(x),v_2(x)]\in L^2(0,\pi,E^2)$.
Domain $D$ consists of the functions $y(x)\in W_1^2(0,\pi,E^2)$ satisfying the boundary conditions $\langle y(0),\Phi(\alpha)\rangle=0$, $\langle y(\pi),\Phi(\beta)\rangle=0$ ($\Phi(\alpha)=\col[\sin\alpha,\cos\alpha]$). Spectrums of such operators are described. Method for solving the inverse spectrum problem is given, viz., technique to reconstruct number $\mu$ and function $v(x)$ from the spectral data is presented. For the case of real $v(x)$, it is shown that the number $\mu$ and the non-local potential $v(x)$ are unambiguously defined by the two spectra. Complete description of spectral data of inverse problem for this class of operators is given.
\vspace{5mm}

\section{Characteristic function. Eigenfunctions}\label{s1}

{\bf 1.1.} In the Hilbert space $L^2(0,\pi,E^2)$ of complex vector-valued functions $f(x)=\col[f_1(x),f_2(x)]\in E^2$ ($x\in[0,\pi]$), consider the free Dirac operator $D_0$,
\begin{equation}
(D_0y)(x)\stackrel{\d}{=}\left[
\begin{array}{ccc}
0&1\\
-1&0
\end{array}\right]\frac{dy(x)}{dx}\quad(x\in(0,\pi)),\label{eq1.1}
\end{equation}
with the domain
\begin{equation}
\mathfrak{D}(D_0)=\{y(x)\in W_1^2(0,\pi,E^2):\langle y(0),\Phi(\alpha)\rangle=0;\langle y(\pi),\Phi(\beta)\rangle=0\}\label{eq1.2}
\end{equation}
where
\begin{equation}
\Phi(\alpha)=\col[\sin\alpha,\cos\alpha]\quad(\alpha,\beta\in\mathbb{R}).\label{eq1.3}
\end{equation}
Eigenfunctions $D_0y=\lambda y$ of the operator $D_0$ satisfy the system of equations
\begin{equation}
\left\{
\begin{array}{lll}
y'_2=\lambda y_1;\\
-y'_1=\lambda y_2.
\end{array}\right.\label{eq1.4}
\end{equation}
It is easy to see that the vector-valued function
\begin{equation}
y_\alpha(\lambda,x)=\col[\cos(\lambda x-\alpha),\sin(\lambda x-\alpha)]\label{eq1.5}
\end{equation}
is the solution to system \eqref{eq1.4} and satisfies the first boundary condition $\langle y_\alpha(\lambda,0),$ $\Phi(\alpha)\rangle=0$. The second boundary condition for $y_\alpha(x,\lambda)$ implies
$$0=\langle y_\alpha(\lambda,\pi),\Phi(\beta)\rangle=\cos(\lambda\pi-\alpha)\sin\beta+\sin(\lambda\pi-\alpha)\cos\beta=\sin(\lambda\pi-\alpha+\beta).$$
The function
\begin{equation}
\Delta(0,\lambda)=\sin(\lambda\pi-\alpha+\beta)\label{eq1.6}
\end{equation}
is said to be the characteristic function of the operator $D_0$ \eqref{eq1.1}, \eqref{eq1.2}. Zeros $\lambda=\lambda_k(0)$ of the function $\Delta(0,\lambda)$ define spectrum of the operator $D_0$,
\begin{equation}
\sigma(D_0)=\left\{\lambda_k(0)=k+\frac{\alpha-\beta}\pi:k\in\mathbb{Z}\right\}.\label{eq1.7}
\end{equation}
Consequently, operator $D_0$ \eqref{eq1.1}, \eqref{eq1.2} has simple discrete spectrum and its complete set of the orthonormal eigenfunctions is
\begin{equation}
u(\lambda_k(0),x)=\frac1{\sqrt\pi}\col[\cos(\lambda_k(0)x-\alpha),\sin(\lambda_k(0)x-\alpha)]\quad(k\in\mathbb{Z}).\label{eq1.8}
\end{equation}
\vspace{5mm}

{\bf 1.2.} Denote by $D$ the Dirac operator with non-local potential,
\begin{equation}
(Dy)(x)=(D(\mu,v)y)(x)\stackrel{\d}{=}(D_0y)(x)+\mu\int\limits_0^\pi\langle y(t),v(t)\rangle_{E^2}dtv(x),\label{eq1.9}
\end{equation}
where $\mu\in\mathbb{R}$ and $v=\col[v_1,v_2]\in L^2(0,\pi,E^2)$. The operator $D$ is a one-dimensional perturbation of the operator $D_0$ \eqref{eq1.1}, \eqref{eq1.2} and $\left. D_\mu\right|_{\mu=0}=D_0$ \eqref{eq1.1}. Domains of $D_0$ and $D$ coincide.

It is easy to see that solution to the equation
$$D_0y-\lambda y=f$$
satisfying the boundary condition $\langle y(0),\Phi(\alpha)=0$ is
\begin{equation}
y(\lambda,x)=a(\lambda)y_\alpha(\lambda,x)+\int\limits_0^x\left[
\begin{array}{ccc}
\cos\lambda(x-t)&-\sin\lambda(x-t)\\
\sin\lambda(x-t)&\cos\lambda(x-t)
\end{array}\right]\left(
\begin{array}{ccc}
-f_2(t)\\
f_1(t)
\end{array}\right)dt\label{eq1.10}
\end{equation}
where $a(\lambda)$ is an arbitrary function of $\lambda$. Hence it follows that the eigenfunction $Du=\lambda u$ of the operator $D$ \eqref{eq1.9}, which is the solution to the equation
$$D_0y+\mu\langle y,v\rangle v=\lambda u$$
and for which $\langle u(0),\Phi(\alpha)\rangle=0$, satisfies the integral equation
\begin{equation}
u(\lambda,x)=a(\lambda)y_\alpha(\lambda,x)-b(\lambda)\int\limits_0^x\left[
\begin{array}{ccc}
\cos\lambda(x-t)&-\sin\lambda(x-t)\\
\sin\lambda(x-t)&\cos\lambda(x-t)
\end{array}\right]\left(
\begin{array}{ccc}
-v_2(t)\\
v_1(t)
\end{array}\right)dt,\label{eq1.11}
\end{equation}
where $b(\lambda)=\mu\langle u(\lambda,x),v(x)\rangle$. Scalarly multiplying this equality by $\mu v(x)$, one obtains
\begin{equation}
\mu a(\lambda)\langle y_\alpha,v\rangle-b(\lambda)(\mu F(\lambda)+1)=0\label{eq1.12}
\end{equation}
where
\begin{equation}
F(\lambda)=\int\limits_0^\pi\left\langle\int\limits_0^x\left[
\begin{array}{ccc}
\cos\lambda(x-t)&-\sin\lambda(x-t)\\
\sin\lambda(x-t)&\cos\lambda(x-t)
\end{array}\right]\left(
\begin{array}{ccc}
-v_2(t)\\
v_1(t)
\end{array}\right)dt,v(x)\right\rangle dx.\label{eq1.13}
\end{equation}
Taking into account that $\langle y_\alpha(\lambda,\pi),\Phi(\beta)\rangle=\Delta(0,\lambda)$, from the second boundary condition $\langle u(\lambda,\pi,\Phi(\beta)\rangle=0$ one obtains the equality
\begin{equation}
a(\lambda)\Delta(0,\lambda)-b(\lambda)G(\lambda)=0,\label{eq1.14}
\end{equation}
where
\begin{equation}
G(\lambda)=\left\langle\int\limits_0^\pi\left[
\begin{array}{ccc}
\cos\lambda(\pi-t)&-\sin\lambda(\pi-t)\\
\sin\lambda(\pi-t)&\cos\lambda(\pi-t)
\end{array}\right]\left(
\begin{array}{ccc}
-v_2(t)\\
v_1(t)
\end{array}\right)dt,\Phi(\beta)\right\rangle.\label{eq1.15}
\end{equation}
Relations \eqref{eq1.12} and \eqref{eq1.14} give the system of linear equations
$$\left\{
\begin{array}{lll}
a(\lambda)\mu\langle y_\alpha,v\rangle-b(\lambda)(\mu F(\lambda)+1)=0;\\
a(\lambda)\Delta(0,\lambda)-b(\lambda)G(\lambda)=0;
\end{array}\right.$$
for $a(\lambda)$, $b(\lambda)$, which has non-trivial solutions if determinant of this system $\Delta(\mu,\lambda)$ vanishes, $\Delta(\mu,\lambda)=0$, where
\begin{equation}
\Delta(\mu,\lambda)=\Delta(0,\lambda)+\mu\{\Delta(0,\lambda)\cdot F(\lambda)-\langle y_\alpha,v\rangle\cdot G(\lambda)\}.\label{eq1.16}
\end{equation}
Function $\Delta(\mu,\lambda)$ is said to be the characteristic function of the operator $D$ \eqref{eq1.9}, its zeros $\lambda_k(\mu)$ form the spectrum of operator $D$,
\begin{equation}
\sigma(D)=\{\lambda_k(\mu):\Delta(\mu,\lambda_k(\mu))=0\}.\label{eq1.17}
\end{equation}
\vspace{5mm}

{\bf 1.3.} Transform expression \eqref{eq1.16} for the characteristic function $\Delta(\mu,\lambda)$. From \eqref{eq1.15} one has
$$G(\lambda)=\left\langle\int\limits_0^\pi\left(
\begin{array}{ccc}
-v_2(t)\cos\lambda(\pi-t)-v_1(t)\sin\lambda(\pi-t)\\
-v_2(t)\sin\lambda(\pi-t)+v_1(t)\cos\lambda(\pi-t)
\end{array}\right)dt\left(
\begin{array}{ccc}
\sin\beta\\
\cos\beta
\end{array}\right)\right\rangle$$
$$=\int\limits_0^\pi\left\{v_1(t)\cos(\lambda(\pi-t)+\beta)-v_2(t)\sin(\lambda(\pi-t)+\beta)\right\}dt$$
$$=\int\limits_0^\pi\{v_1(t)\cos(\lambda(t-\pi)-\beta)+v_2(t)\sin(\lambda(t-\pi)-\beta)\}dt$$
$$=\frac12\int\limits_0^\pi\left\{e^{-i(\lambda\pi+\beta)}\int\limits_0^\pi e^{i\lambda t}(v_1(t)-iv_2(t))dt\right.$$
$$\left.+e^{i(\lambda\pi+\beta)}\int\limits_0^\pi e^{-i\lambda t}(v_1(t)+iv_2(t))dt\right\}.$$
In terms of Fourier transforms
\begin{equation}
\widetilde{v}_\pm(\lambda)\stackrel{\d}{=}\int\limits_0^\pi e^{-i\lambda x}v_\pm(x)dx\label{eq1.18}
\end{equation}
of the functions
\begin{equation}
v_\pm(x)=v_1(x)\pm iv_2(x),\label{eq1.19}
\end{equation}
expression for $G(\lambda)$ is
\begin{equation}
G(\lambda)=\frac12\left\{e^{i(\lambda\pi+\beta)}\widetilde{v}_+(\lambda)+e^{-i(\lambda\pi+\beta)}\widetilde{v}_-(-\lambda)\right\}.\label{eq1.20}
\end{equation}
Since
$$\langle y_\alpha,v\rangle=\int\limits_0^\pi[\cos(\lambda x-\alpha)\overline{v_1}(x)+\sin(\lambda x-\alpha)\overline{v_2}(x)]dx$$
$$=\frac12\int\limits_0^\pi\left\{e^{-i\alpha}\int\limits_0^\pi e^{i\lambda x}(\overline{v_1}(x)-i\overline{v_2}(x))dx+e^{ix}\int\limits_0^\pi e^{-i\lambda x}(\overline{v_1}(x)+i\overline{v_2(x)})dx\right\}$$
and according to \eqref{eq1.18}, \eqref{eq1.19}
\begin{equation}
\langle y_\alpha,v\rangle=\frac12\left\{e^{-i\alpha}\widetilde{v}_+^*(\lambda)+e^{i\alpha}\widetilde{v}_-^*(-\lambda)\right\}\label{eq1.21}
\end{equation}
where
$$f^*(\lambda)=\overline{f(\overline\lambda)}.$$
Relation \eqref{eq1.13} implies
$$F(\lambda)=\int\limits_0^\pi\left\langle\int\limits_0^x\left(
\begin{array}{ccc}
-v_2(t)\cos\lambda(x-t)-v_1(t)\sin\lambda(x-t)\\
-v_2(t)\sin\lambda(x-t)+v_1(t)\cos\lambda(x-t)
\end{array}\right)dt,\left(
\begin{array}{ccc}
v_1(x)\\
v_2(x)
\end{array}\right)\right\rangle dx$$
$$=\int\limits_0^\pi\int\limits_0^x\{\cos\lambda(x-t)[v_1(t)\overline{v_2}(x)-v_2(t)\overline{v_1(x)}]$$
$$+\sin\lambda(x-t)[-v_1(t)\overline{v_1(x)}-v_2(t)\overline{v_2(x)}]\}dtdx$$
$$=\int\limits_0^\pi dx\int\limits_0^xdt\left\{\cos\lambda(x-t)\left\langle\left[
\begin{array}{ccc}
0&-1\\
1&0
\end{array}\right]v(t),v(x)\right\rangle-\sin\lambda(x-t)\langle v(t),v(x)\rangle\right\}$$
$$=\frac12\int\limits_0^\pi dx\int\limits_0^xdt\left\{e^{i\lambda(x-t)}\left\langle\left[
\begin{array}{ccc}
i&-1\\
1&i
\end{array}\right]v(t),v(x)\right\rangle\right.$$
$$\left.+e^{-i\lambda(x-t)}\left\langle\left[
\begin{array}{ccc}
-i&-1\\
1&-i
\end{array}\right]v(t),v(x)\right\rangle\right\}.$$
Taking into account
$$\left\langle\left[
\begin{array}{ccc}
i&-1\\
1&i
\end{array}\right]v,w\right\rangle=\left\langle\left(
\begin{array}{ccc}
i(v_1+iv_2)\\
v_1+iv_2
\end{array}\right),w\right\rangle=(v_1+iv_2)(i\overline w_1+\overline w_2)$$
$$=i(v_1+iv_2)(\overline{w_1+iw_2}),$$
one has
$$F(\lambda)=\frac i2\left\{\int\limits_0^\pi dx\int\limits_0^xdte^{i\lambda(x-t)}v_+(t)\overline v_+(x)-\int\limits_0^\pi dx\int\limits_0^xdte^{-i\lambda(x-t)+}v_-(t)\overline{v_-(x)}\right\}.$$
Define the functions
\begin{equation}
\Phi_\pm(\lambda)=\int\limits_0^\pi dx\int\limits_0^xe^{-i\lambda(x-t)}\overline{v_\pm(t)}dtv_\pm(x),\label{eq1.22}
\end{equation}
then
\begin{equation}
F(\lambda)=\frac i2\{\Phi_+^*(\lambda)-\Phi_-^*(-\lambda)\}.\label{eq1.23}
\end{equation}
Notice that
\begin{equation}
\Phi_\pm(\lambda)=\int\limits_0^\pi dx\int\limits_0^xe^{-i\lambda\xi}\overline{v_\pm(x-\xi)}d\xi v_\pm(x)=\int\limits_0^\pi d\xi e^{-i\lambda\xi}g_\pm(\xi),\label{eq1.24}
\end{equation}
where
\begin{equation}
g_\pm(\xi)=\int\limits_\xi^\pi\overline{v_\pm(x-\xi)}v_\pm(x)dx.\label{eq1.25}
\end{equation}
So, $\Phi_\pm(\lambda)$ is the Fourier transform \eqref{eq1.24} of the convolution $g_\pm$ \eqref{eq1.25}. Since $v_\pm\in L^2(0,\pi)$ \cite{4}, $g_\pm\in L^2(0,\pi)$. Therefore $\widetilde v_\pm(\lambda)$ and $\Phi_\pm(\lambda)$ are entire functions of exponential type ($\sigma\leq\pi$) due to the Paley -- Wiener theorem \cite{4}. The following statement is an analogue of the well-known statement on the fact that Fourier transform of a convolution from $L^2(\mathbb{R})$ equals the product of Fourier transforms of elements of the convolution.

\begin{lemma}\label{l1.1}
For functions $\Phi_\pm(\lambda)$ \eqref{eq1.22}, the identity
\begin{equation}
\Phi_\pm(\lambda)+\Phi_\pm^*(\lambda)=\widetilde v_\pm(\lambda)\widetilde v_\pm^*(\lambda)\label{eq1.26}
\end{equation}
holds, where $\widetilde v_{\pm}(\lambda)$ are given by \eqref{eq1.19}.
\end{lemma}

P r o o f. Equality \eqref{eq1.26} follows from the integration by parts formula,
$$\Phi_\pm(\lambda)=\int\limits_0^\pi\int\limits_0^xe^{i\lambda t}\overline{v_\pm(t)}dtd\int\limits_0^xe^{-i\lambda s}v_\pm(s)ds=\widetilde{v}_\pm(\lambda)\widetilde v_\pm^*(\lambda)$$
$$-\int\limits_0^\pi\int\limits_0^xe^{-i\lambda s}v_\pm(s)dse^{i\lambda x}\overline{v_\pm(x)}dx=\widetilde v_\pm(\lambda)\widetilde v_\pm^*(\lambda)-\Phi_\pm^*(\lambda). \quad \blacksquare$$

\begin{remark}
It is easy to show that the entire function $\Phi_\pm(\lambda)$ of exponential type ($\sigma\leq\pi$) satisfying \eqref{eq1.26} is the Fourier transform \eqref{eq1.24} of the convolution $g_\pm$ \eqref{eq1.25} where $\widetilde v_\pm(\lambda)$ are given by \eqref{eq1.19}.
\end{remark}

Substituting expressions for $G(\lambda)$ \eqref{eq1.20}, $\langle y_\alpha,v\rangle$ \eqref{eq1.21} and $F(\lambda)$ \eqref{eq1.23} into \eqref{eq1.16}, one obtains
$$\Delta(\mu,\lambda)=\Delta(0,\lambda)+\frac\mu 4\left\{\left(e^{i(\lambda\pi+\beta-\alpha)}-e^{-i(\lambda\pi+\beta-\alpha)}\right)(\Phi_+^*(\lambda)-\Phi_-^*(-\lambda))\right.$$
$$\left.-\left(e^{-i\alpha}\widetilde v_+^*(\lambda)+e^{i\alpha}\widetilde v_-^*(-\lambda)\right)\left(e^{i(\lambda\pi+\beta)}\widetilde v_+(\lambda)+e^{-i(\lambda\pi+\beta)}\widetilde v_-(-\lambda)\right)\right\}$$
$$=\Delta(0,\lambda)+\frac\mu 4\left\{\left(e^{i(\lambda\pi+\beta-\alpha)}-e^{-i(\lambda\pi+\beta-\alpha)}\right)\Phi_+^*(\lambda)-\left(e^{i(\lambda\pi+\beta-\alpha)}-e^{-i(\lambda\pi+\beta-\alpha)}\right)\right.$$
$$\times\Phi_-^*(-\lambda)-e^{i(\lambda\pi+\beta-\alpha)}\widetilde v_+(\lambda)\widetilde v_+^*(\lambda)-e^{-i(\lambda\pi+\beta-\alpha)}\widetilde v_-(-\lambda)\widetilde v_-(-\lambda)$$
$$\left.-e^{i(\lambda\pi+\beta+\alpha)}\widetilde v_+(\lambda)\widetilde v_-^*(-\lambda)-e^{-i(\lambda\pi+\beta+\alpha)}\widetilde v(-\lambda)\widetilde v_+^*(\lambda)\right\},$$
using \eqref{eq1.26} one obtains
$$\Delta(\mu,\lambda)=\Delta(0,\lambda)-\frac\mu 4\left\{e^{i(\lambda\pi+\beta-\alpha)}\Phi_+(\lambda)+e^{-i(\lambda\pi+\beta-\alpha)}\Phi_+^*(\lambda)\right.$$
$$+e^{i(\lambda\pi+\beta-\alpha)}\Phi_-^*(-\lambda)+e^{-i(\lambda\pi+\beta-\alpha)}\Phi_-(-\lambda)+e^{i(\lambda\pi+\beta+\alpha)}\widetilde v_+(\lambda)\widetilde v_-^*(-\lambda)$$
$$\left.+e^{-i(\lambda\pi+\beta+\alpha)}\widetilde v_-(-\lambda)\widetilde v_+^*(\lambda)\right\}=\Delta(0,\lambda)-\frac\mu 4\left\{e^{i(\lambda\pi+\beta-\alpha)}\right.$$
$$\times\left[\Phi_+(\lambda)+\Phi_-^*(-\lambda)+e^{2i\alpha}\widetilde v_+(\lambda)\widetilde v_-^*(-\lambda)\right]$$
$$\left.+e^{-i(\lambda\pi+\beta-\alpha)}\left[\Phi_+^*(\lambda)+\Phi_-(-\lambda)+e^{-2i\alpha}\widetilde v_-(-\lambda)\widetilde v_+^*(\lambda)\right]\right\}.$$

\begin{theorem}\label{t1.1}
Characteristic function $\Delta(\mu,\lambda)$ \eqref{eq1.16} is expressed via the functions $\widetilde v_\pm(\lambda)$ \eqref{eq1.19} and $\Phi_\pm(\lambda)$ \eqref{eq1.22} by the formula
\begin{equation}
\Delta(\mu,\lambda)=\Delta(0,\lambda)-\frac\mu 4\left\{e^{i(\lambda\pi+\beta-\alpha)}R(\lambda)+e^{-i(\lambda\pi+\beta-\alpha)}R^*(\lambda)\right\}\label{eq1.27}
\end{equation}
where
\begin{equation}
R(\lambda)=\Phi_+(\lambda)+\Phi_-^*(-\lambda)+e^{2i\alpha}\widetilde v_+(\lambda)\widetilde v_-^*(-\lambda).\label{eq1.28}
\end{equation}
\end{theorem}

Relation \eqref{eq1.27} implies
\begin{equation}
\Delta^*(\mu,\lambda)=\Delta(\mu,\lambda).\label{eq1.29}
\end{equation}
\vspace{5mm}

{\bf IV.} Find the form of eigenfunctions of the operator $D$ \eqref{eq1.9}. If $\lambda=\lambda_k(\mu)$ is a root of $\Delta(\mu,\lambda)$ \eqref{eq1.27}, then from \eqref{eq1.14} one finds
$$a(\lambda)=\frac{b(\lambda)}{\Delta(0,\lambda)}G(\lambda).$$
Substituting this expression into \eqref{eq1.11}, one finds eigenfunctions of the operator $D$.

\begin{theorem}\label{t1.2}
To each zero $\lambda=\lambda_k(\mu)$ of the characteristic function $\Delta(\mu,\lambda)$ \eqref{eq1.27} there corresponds an eigenfunction $u(\lambda,x)$ of the operator $D$ \eqref{eq1.9},
\begin{equation}
\begin{array}{ccc}
u(\lambda,x)=-\left(e^{i(\lambda\pi+\beta)}\widetilde v_+(\lambda)+e^{-i(\lambda\pi+\beta)}\widetilde v_-(-\lambda)\right)y_\alpha(\lambda,x)\\
{\displaystyle-\sin(\lambda\pi-\alpha+\beta)\left\{\left(
\begin{array}{ccc}
i\\
1
\end{array}\right)\int\limits_0^xe^{i\lambda(x-t)}v_+(t)dt+\left(
\begin{array}{ccc}
-i\\
1
\end{array}\right)\int\limits_0^xe^{-i\lambda(x-t)}v_-(t)dt\right\}}
\end{array}\label{eq1.30}
\end{equation}
where $v_\pm(t)$ and $\widetilde v_\pm(\lambda)$ are given by \eqref{eq1.18} and \eqref{eq1.19} correspondingly and $y_\alpha(\lambda,x)$ is given by \eqref{eq1.5}.
\end{theorem}

P r o o f. Representation \eqref{eq1.30} follows from \eqref{eq1.11} due to \eqref{eq1.20}. The second term, up to $\Delta(0,\lambda)$, equals
$$\int\limits_0^x\left[
\begin{array}{ccc}
\cos\lambda(x-t)&-\sin\lambda(x-t)\\
\sin\lambda(x-t)&\cos\lambda(x-t)
\end{array}
\right]\left(
\begin{array}{ccc}
-v_2(t)\\
v_1(t)
\end{array}
\right)dt$$
$$=\frac12\int\limits_0^x\left[
\begin{array}{ccc}
e^{i\lambda(x-t)}+e^{-i\lambda(x-t)}&i\left(e^{i\lambda(x-t)}-e^{-i\lambda(x-t)}\right)\\
-i\left(e^{i\lambda(x-t)}-e^{-i\lambda(x-t)}\right)&e^{i\lambda(x-t)}+e^{i\lambda(x-t)}
\end{array}\right]\left(
\begin{array}{ccc}
-v_2(t)\\
v_1(t)
\end{array}\right)dt$$
$$=\frac12\int\limits_0^x\left\{e^{i\lambda(x-t)}\left[
\begin{array}{ccc}
1&i\\
-i&1
\end{array}\right]+e^{-i\lambda(x-t)}\left[
\begin{array}{ccc}
1&-i\\
i&1
\end{array}
\right]\right\}\left(
\begin{array}{ccc}
-v_2(t)\\
v_1(t)
\end{array}\right)dt$$
$$=\frac12\int\limits_0^x\left\{
e^{i\lambda(x-t)}\left(
\begin{array}{ccc}
i\\
1
\end{array}\right)v_+(t)+e^{-i\lambda(x-t)}\left(
\begin{array}{ccc}
-i\\
1
\end{array}\right)v_-(t)\right\}dt$$
this proves \eqref{eq1.30}.

Note that $b(\lambda)=0$, when $\lambda_k(\mu)=\lambda_s(0)$, and the function $u(\lambda,x)$ \eqref{eq1.30} is proportional to the eigenfunction $u(\lambda_k(0),x)$ \eqref{eq1.8} of the operator $D_0$.
\vspace{5mm}

\section{Spectrum of operator $D$}\label{s2}

{\bf 2.1.} Denote by $L_0$ the selfadjoint operator with dense domain $\mathfrak{D}(L_0)$ given in a Hilbert space $H$, which has simple discrete spectrum
\begin{equation}
\sigma(L_0)=\{\lambda_k:k\in\mathbb{Z}\}\label{eq2.1}
\end{equation}
assuming that real numbers $\lambda_k$ do not have finite limit points and are enumerated in ascending order. And let $u_k$ be orthonormal $\langle u_k,u_s\rangle=\delta_{k,s}$ eigenfunctions $L_0u_k=\lambda_ku_k$, of the operator $L_0$. Spectrum expansion of $L_0$ is
$$L_0=\sum\limits_k\lambda_k\langle\cdot,u_k\rangle u_k,$$
and its resolvent $R_{L_0}(\lambda)=(L_0-\lambda I)^{-1}$ equals
\begin{equation}
R_{L_0}(\lambda)f=\sum\limits_k\frac{f_k}{\lambda_k-\lambda}u_k,\label{eq2.2}
\end{equation}
where $f_k=\langle f,u_k\rangle$.

Consider the selfadjoint operator (cf. \eqref{eq1.9}),
\begin{equation}
L=L_0+\mu\langle\cdot,v\rangle v,\label{eq2.3}
\end{equation}
where $\mu\in\mathbb{R}$ and $v\in H$, which is a one-dimensional perturbation of $L_0$. Domains of $L$ and $L_0$ coincide, $\mathfrak{D}(L)=\mathfrak{D}(L_0)$. Resolvent $R_L(\lambda)=(L-\lambda I)^{-1}$ of the operator $L$ \eqref{eq2.3} is expressed via the resolvent $R_{L_0}(\lambda)$ by the formula
\begin{equation}
R_L(\lambda)f=R_{L_0}(\lambda)f-\frac{\mu\langle R_{L_0}f,v\rangle}{1+\mu\langle R_{L_0}(\lambda)v,v\rangle}R_{L_0}(\lambda)v\label{eq2.4}
\end{equation}
or, due to \eqref{eq2.2},
\begin{equation}
R_L(\lambda)f=\sum\limits_k\left\{\frac{f_k}{\lambda_k-\lambda}-\frac{\displaystyle\mu\sum\limits_s\frac{f_s\overline{u_s}}{\lambda_s-\lambda}}{\displaystyle1+\mu\sum\limits_s\frac{|v_s|^2}{\lambda_s-\lambda}}\cdot\frac{v_k}{
\lambda_k-\lambda}\right\}u_k.\label{eq2.5}
\end{equation}
Hence it follows that if $v_p=0$, then the point $\lambda_p$ belongs to the spectrum of operator $L$. If $v_p\not=0$, then calculating residue of $R_L(\lambda)f$ at the point $\lambda=\lambda_p$, one obtains
$$\res\limits_{\lambda_p}R_L(\lambda)f=f_p-\frac{\mu f_p\overline{v_p}}{\mu|v_p|^2}v_p=0,$$
and thus $R_L(\lambda)f$ has no singularity at the point $\lambda=\lambda_p$. Therefore it is natural to divide the set $\sigma(L_0)$ \eqref{eq2.1} into two disjoint subsets, $\sigma(L_0)=\sigma_0\cup\sigma_1$ ($\sigma_0\cap\sigma_1=\emptyset$),
\begin{equation}
\sigma_0\stackrel{\d}{=}\{\lambda_k\in\sigma(L_0):v_k=0\};\quad\sigma_1\stackrel{\d}{=}\{\lambda_k\in\sigma(L_0):v_k\not=0\}.\label{eq2.6}
\end{equation}
Singularities of the resolvent $R_L(\lambda)$ \eqref{eq2.5} contain zeros of the function $Q(\lambda)$,
\begin{equation}
Q(\lambda)\stackrel{\d}{=}1+\mu\langle R_{L_0}(\lambda)v,v\rangle=1+\mu\sum\limits_s\frac{|v_s|^2}{\lambda_s-\lambda},\label{eq2.7}
\end{equation}
besides, summation is carried out by those $s$ for which $\lambda_s\in\sigma_1$.

\begin{remark}\label{r2.1}
Zeros of $Q(\lambda)$ are real, simple and alternate with points from $\sigma_1$, besides, if $\mu>0$ ($<0$), then zeros of $Q(\lambda)$ translate to the right (left) compared with the points $\lambda_k\in\sigma_1$.
\end{remark}

Denote by $\sigma_2$ the set of zeros of $Q(\lambda)$ \eqref{eq2.7},
\begin{equation}
\sigma_2\stackrel{\d}{=}\{\mu_k:Q(\mu_k)=0\}.\label{eq2.8}
\end{equation}
The numbers $\mu_k$ can coincide with the elements $\lambda_s$ from $\sigma_0$.

\begin{theorem}\label{t2.1}
Let $L_0$ be a selfadjoint operator with simple discrete spectrum $\sigma(L_0)$ \eqref{eq2.1} and $\sigma(L_0)$ do not have finite limit points. Then spectrum of the operator $L$ \eqref{eq2.3} equals
\begin{equation}
\sigma(L)=\{\sigma_0\backslash(\sigma_0\cap\sigma_2)\}\cup\{\sigma_2\backslash(\sigma_0\cap\sigma_2)\}\cup\{\sigma_0\cap\sigma_2\},\label{eq2.9}
\end{equation}
where points from $\sigma_0\backslash(\sigma_0\cap\sigma_2)$, $\sigma_2\backslash(\sigma_0\cap\sigma_2)$ are of multiplicity $1$ and the points from $\sigma_0\cap\sigma_2$ are of multiplicity $2$.
\end{theorem}

\begin{remark}\label{r2.2}
If spectrum $\sigma(L_0)$ \eqref{eq2.1} has {\bf separability property}, i. e.,
\begin{equation}
d=\inf\limits_k(\lambda_{k+1}-\lambda_k)>0,\label{eq2.10}
\end{equation}
then the set $\sigma_0\cap\sigma_2$ is finite.
\end{remark}

If this is not the case, then there exists the sequence of zeros $\mu_s$ of the function $Q(\lambda)$ which coincides with $\lambda_p\in\sigma_0$, and
\begin{equation}
1+\mu\sum\limits_s\frac{|v_s|^2}{\lambda_s-\lambda_p}=0\quad(\lambda_p=\mu_s\in\sigma_0\cap\sigma_2)\label{eq2.11}
\end{equation}
where summation is carried out by $\lambda_s\in\sigma_1$ ($\sigma_1\cap\sigma_0=\emptyset$). Since $|\lambda_s-\lambda_p|\geq d$ and $\sum\limits_s|v_s|^2<\infty$ ($v\in H$), series \eqref{eq2.11} converges uniformly by $p$. As $p$ approaches infinity, one obtains the contradiction, $1=0$.

Note that the operator $D_0$ \eqref{eq1.1}, \eqref{eq1.2} has the separability property \eqref{eq2.10} due to the form of $\sigma(D_0)$ \eqref{eq1.7}.
\vspace{5mm}

{\bf 2.2.} Calculate the function $Q(\lambda)$ \eqref{eq2.7} for the operator $L_0=D_0$ \eqref{eq1.1}, \eqref{eq1.2}. Using \eqref{eq1.10} and method of variation of arbitrary constants $(a(\lambda))$, one finds the resolvent $R_{D_0}(\lambda)=(D_0-\lambda I)^{-1}$ of the operator $D_0$ \cite{1,2},
\begin{equation}
\begin{array}{ccc}
{\displaystyle(R_{D_0}(\lambda)f)(x)=-\frac1{\Delta(0,\lambda)}\left\{y_\beta(\lambda,x-\pi)\int\limits_0^x\langle f(t),y_\alpha(\overline\lambda,t)\rangle_{E^2}dt\right.}\\
\left.+y_\alpha(\lambda,x)\int\limits_x^\pi\langle f(t),y_\beta(\overline\lambda,t-\pi)\rangle_{E^2}dt\right\}
\end{array}\label{eq2.12}
\end{equation}
where $y_\alpha(\lambda,x)$ and $\Delta(0,\lambda)$ are given by \eqref{eq1.5} and \eqref{eq1.6}. Hence it follows that
\begin{equation}
\begin{array}{ccc}
{\displaystyle\langle R_{D_0}(\lambda)v,v\rangle=-\frac1{\Delta(0,\lambda)}\left\{\int\limits_0^\pi\int\limits_0^x\langle v(t),y_\alpha(\overline\lambda,t)\rangle dt\langle y_\beta(\lambda,x-\pi),v(x)\rangle dx\right.}\\
{\displaystyle\left.+\int\limits_0^\pi\int\limits_x^\pi\langle v(t),y_\beta(\overline\lambda,t-\pi)\rangle dt\langle y_\alpha(\lambda,x),v(x)\rangle dx\right\},}
\end{array}\label{eq2.13}
\end{equation}
and since
$$\langle v(t),y_\alpha(\overline\lambda,t)=v_1(t)\cos(\lambda t-\alpha)+v_2(t)\sin\lambda(t-\alpha)$$
\begin{equation}
\begin{array}{ccc}
{\displaystyle=\frac12\left\{e^{i(\lambda t-\alpha)}v_-(t)+e^{-i(\lambda t-\alpha)}v_+(t)\right\};}\\
\langle v(t),y_\beta(\overline\lambda,t-\pi)\rangle=v_1(t)\cos(\lambda(t-\pi)-\beta)+v_2(t)\sin(\lambda(t-\pi)-\beta)
\end{array}\label{eq2.14}
\end{equation}
$$=\frac12\left\{e^{i(\lambda t-\lambda\pi-\beta)}v_-(t)+e^{-i(\lambda t-\lambda\pi-\beta)}v_+(t)\right\},$$
then
$$\langle R_{D_0}(\lambda)v,v\rangle=-\frac1{4\cdot\Delta(0,\lambda)}\cdot\left\{\int\limits_0^\pi dx\int\limits_0^xdt\left(e^{i(\lambda t-\alpha)}v_-(t)+e^{-i(\lambda t-\alpha)}v_+(t)\right)\right.$$
$$\times\left(e^{-i(\lambda x-\lambda\pi-\beta)}\overline{v_-(x)}+e^{i(\lambda x-\lambda\pi-\beta)}\overline{v_+(x)}+\int\limits_0^\pi dx\int\limits_x^\pi dt\right.$$
$$\left.\times\left(e^{i(\lambda t-\lambda\pi-\beta)}v_-(t)+e^{-i(\lambda t-\lambda\pi-\beta)}v_+(t)\right)\left(e^{-i(\lambda x-\alpha)}\overline{v_-(x)}+e^{i(\lambda x-\alpha)}\overline{v_+(x)}\right)\right\}$$
$$=-\frac1{4\Delta(0,\lambda)}\int\limits_0^\pi dx\left\{\int\limits_0^xdt\left[e^{i\lambda(t-x)}e^{i(\lambda\pi+\beta-\alpha)}v_-(t)\overline{v_-(x)}\right.\right.$$
$$+e^{i\lambda(x-t)}e^{-i(\lambda\pi+\beta-\alpha)}v_+(t)\overline{v_+(x)}+e^{i\lambda(t+x)}e^{-i(\lambda\pi+\beta+\alpha)}v_-(t)\overline{v_+(x)}$$
$$\left.+e^{-i\lambda(t+x)}e^{i(\lambda\pi+\beta+\alpha)}v_+(t)\overline{v_-(x)}\right]+\int\limits_x^\pi dt\left[e^{i\lambda(t-x)}e^{-i(\lambda\pi+\beta-\alpha)}v_-(t)\overline{v_-(x)}\right.$$
$$+e^{i\lambda(x-t)}e^{i(\lambda\pi+\beta-\alpha)}v_+(t)\overline{v_+(x)}+e^{i\lambda(t+x)}e^{-i(\lambda\pi+\beta+\alpha)}v_-(t)\overline{v_+(x)}$$
$$\left.+e^{-i(x+t)}e^{i(\lambda\pi+\beta+\alpha)}v_+(t)\overline{v_-(x)}\right].$$
Using \eqref{eq1.19} and \eqref{eq1.22}, one obtains
$$\langle R_{D_0}(\lambda)v,v\rangle=-\frac1{4\Delta(0,\lambda)}\left\{e^{-i(\lambda\pi+\beta-\alpha)}\Phi_+^*(\lambda)+e^{i(\lambda\pi+\beta-\alpha)}\Phi_+(\lambda)\right.$$
$$+e^{i(\lambda\pi+\beta-\alpha)}\Phi_-^*(-\lambda)+e^{-i(\lambda\pi+\beta-\alpha)}\Phi_-(-\lambda)+e^{i(\lambda\pi+\beta+\alpha)}\widetilde v_+(\lambda)\widetilde v_-^*(-\lambda)$$
$$\left.+e^{-i(\lambda\pi+\alpha+\beta)}\widetilde v_-(-\lambda)\widetilde v_+^*(\lambda)\right\}=-\frac1{4\Delta(0,\lambda)}\left\{e^{i(\lambda\pi+\beta-\alpha)}[\Phi_+(\lambda)+\Phi_-^*(-\lambda)\right.$$
$$+e^{2i\alpha}\widetilde{v}_+(\lambda)\widetilde{v}^*(-\lambda)]+e^{-i(\lambda\pi+\beta-\alpha)}[\Phi_+^*(\lambda)+\Phi_-^*(-\lambda)+e^{-2i\alpha}\widetilde{v}_+^*(\lambda)\widetilde{v}_-(-\lambda)]\}$$
$$=-\frac1{4\Delta(0,\lambda)}\{e^{i(\lambda\pi+\beta-\alpha)}R(\lambda)+e^{-i(\lambda\pi+\beta-\alpha)}R^*(\lambda)\},$$
due to \eqref{eq1.28}.

\begin{theorem}\label{t2.2}
The function $Q(\lambda)=1+\mu\langle R_{D_0}(\lambda)v,v\rangle$ \eqref{eq2.7} is expressed via characteristic functions $\Delta(\mu,\lambda)$ \eqref{eq1.27} and
$\Delta(0,\lambda)$ \eqref{eq1.6} by the formula
\begin{equation}\label{eq2.15}
1+\mu\langle R_{D_0}(\lambda)v,v\rangle=\frac{\Delta(\mu,\lambda)}{\Delta(0,\lambda)}.
\end{equation}
\end{theorem}

\begin{remark}\label{r2.3}
Theorem \ref{t2.1} implies that in the quotient $\Delta(\mu,\lambda)/\Delta(0,\lambda)$ zeros corresponding to $\lambda_k(0)\in\sigma_0$ \eqref{eq2.6} truncate and the expression $\Delta(\mu,\lambda)/\Delta(0,\lambda)$ is a meromorphic function, poles of which coincide with the set $\sigma_1$ \eqref{eq2.6} and zeros coincide with the set $\sigma_2$ \eqref{eq2.8}.
\end{remark}

\section{Properties of the characteristic function $\Delta(\mu,\lambda)$}\label{s3}

{\bf 3.1.} A function $f(\lambda)$ is of $C$ (Cartwright) class \cite{5, 6} if

(a) $f(\lambda)$ is an entire function of exponential type;

(b) the integral
\begin{equation}
\int\limits_\mathbb{R}\frac{\ln^+|f(x)|}{1+x^2}<\infty\label{eq3.1}
\end{equation}
is finite. For a function $\Delta(\mu,\lambda)$ \eqref{eq1.27}, the property (a) follows from a Paley -- Wiener theorem and convergence of the integral \eqref{eq3.1} follows from the boundedness of $\Delta(\mu,\lambda)$ as $\lambda\in\mathbb{R}$, therefore the function $\Delta(\mu,\lambda)$ is of $C$ class.

Indicator diagram \cite{5, 6} of the function $\Delta(\mu,\lambda)$ coincides with the segment ${\displaystyle\left[-ih_\Delta\left(-\frac\pi 2\right),ih_\Delta\left(\frac\pi 2\right)\right]}$ of the imaginary axes. Besides, ${\displaystyle h_\Delta\left(-\frac\pi 2\right)=h_\Delta\left(\frac\pi 2\right)}$, this follows from property \eqref{eq1.29} since
$$|\Delta(\mu,iy)|=|\Delta^*(\mu,iy)|=|\overline{\Delta(\mu,\overline{iy})}|=|\Delta(\mu,-iy)|.$$

The following statement \cite{5, 6} is a corollary of the Cartwright -- Levinson theorem.

\begin{theorem}\label{t3.1}
For functions $f(\lambda)$ of $C$ class such that $f(0)\not=0$ and ${\displaystyle h_f\left(\frac\pi 2\right)=}$ ${\displaystyle h_f\left(-\frac\pi 2\right)}$, the following multiplicative representation is true,
$$f(\lambda)=f(0){\rm v.p.}\prod\limits_k\left(1-\frac\lambda{a_k}\right)=f(0)\lim\limits_{R\rightarrow\infty}\prod\limits_{|a_k|<R}\left(1-\frac\lambda{a_k}\right)$$
where $\{a_k\}$ are zeros of the function $f(\lambda)$.
\end{theorem}

Using \eqref{eq1.26}, rewrite the function $\Delta(\mu,\lambda)$ \eqref{eq1.27} as
$$\Delta(\mu,\lambda)=\sin(\lambda\pi-\alpha+\beta)-\frac\mu4\{e^{i(\lambda\pi-\alpha+\beta)}[\Phi_+(\lambda)+\widetilde{v}_-(-\lambda)\widetilde{v}_-^*(-\lambda)$$
$$-\Phi_-(-\lambda)+e^{2i\alpha}\widetilde{v}_+(\lambda)\widetilde{v}_-^*(-\lambda)]+e^{-i(\lambda\pi-\alpha+\beta)}[\widetilde{v}_+(\lambda)\widetilde{v}_+^*(\lambda)-\Phi_+(\lambda)$$
$$+\Phi_-(-\lambda)+e^{-2i\alpha}\widetilde{v}_-(-\lambda)\widetilde{v}_+^*(\lambda)]\},$$
and after elementary transforms one obtains
\begin{equation}
\begin{array}{ccc}
{\displaystyle\Delta(\mu,\lambda)=\sin(\lambda\pi-\alpha+\beta)\left\{1+\frac\mu{2i}(\Phi_+(\lambda)-\Phi_-(-\lambda))\right\}}\\
{\displaystyle-\frac\mu4[\widetilde{v}_+e^{i\alpha}+\widetilde{v}_-(-\lambda)e^{-i\alpha}][\widetilde{v}_+(\lambda)e^{i(\lambda\pi+\beta)}+\widetilde{v}_-(-\lambda)e^{-i(\lambda\pi+\beta)}]^*.}
\end{array}
\end{equation}
Hence it follows that $\Delta(\mu,0)\not=0$, except for rare cases (for example, when $\alpha=\beta$, then $\Delta(\mu,0)=0$ if $\cos\alpha v_1(x)=i\sin\alpha v_2(x)$). So, due to Theorem \ref{t3.1},
$$\Delta(\mu,\lambda)=A{\rm v.p.}\prod\limits_k\left(1-\frac\lambda{\lambda_k(\mu)}\right)\quad(A-{\rm const})$$
and similarly for $\Delta(\mu,0)=\sin(\lambda\pi-\alpha+\beta)$,
$$\Delta(0,\lambda)=\sin(\beta-\alpha)\prod\limits_k\left(1-\frac\lambda{\lambda_k(0)}\right)$$
where $\lambda_k(0)$ are given by \eqref{eq1.7}. Therefore, according to \eqref{eq2.15},
\begin{equation}
1+\mu\sum\limits_k\frac{|v_k|^2}{\lambda_k(0)-\lambda}=\frac A{\sin(\beta-\alpha)}{\rm v.p.}\prod\limits_k\frac{\displaystyle\left(1-\frac \lambda{\lambda_k(\mu)}\right)}{\displaystyle\left(1-\frac\lambda{\lambda_k(0)}\right)},\label{eq3.3}
\end{equation}
besides, $v_k=\langle v(x),u(\lambda_k(0),x)\rangle$ are Fourier coefficients of the vector function $v(x)$ in the basis $u(\lambda_k(0),x)$ \eqref{eq1.8}. Due to Remark \ref{r2.2}, product in \eqref{eq3.3} is realized only by those $k$ for which $\lambda_k\in\sigma_1$ \eqref{eq2.6}.

From equality \eqref{eq3.3} one can find the number $A$,
\begin{equation}
A=\sin(\beta-\alpha)\lim\limits_{y\rightarrow\infty}{\rm v. p.}\prod\limits_k\frac{\displaystyle1-\frac{iy}{\lambda_k(0)}}{\displaystyle1-\frac{iy}{\lambda_k(\mu)}}.\label{eq3.4}
\end{equation}
Calculating residue at the point $\lambda=\lambda_p(0)$ in both parts of equality \eqref{eq3.3}, one obtains
\begin{equation}
\mu|v_p|^2=\frac a{\sin(\beta-\alpha)}\cdot\frac{\lambda_p(0)}{\lambda_p(\mu)}(\lambda_p(\mu)-\lambda_p(0)){\rm v. p.}\prod\limits_{k\not=p}\frac{\lambda_k(0)}{\lambda_k(\mu)}\left(1+\frac{\lambda_k(\mu)-\lambda_p(0)}{\lambda_k(0)-\lambda_p(0)}\right).\label{eq3.5}
\end{equation}

\begin{theorem}\label{t3.2}
The numbers $\mu|v_k|^2$ \eqref{eq3.5} are found unambiguously from the spectrum $\sigma(D)=\{\lambda_k(\mu):k\in\mathbb{Z}\}$ of the operator $D$ \eqref{eq1.9}. Here, $v_k=\langle v(x),u(\lambda_k(0),x)\rangle$ are Fourier coefficients of the function $v(x)$ in the basis of eigenfunctions $u(\lambda_k(0),x)$ \eqref{eq1.8} of the operator $D_0$ \eqref{eq1.1}, \eqref{eq1.2}, besides, the number $A$ is calculated by formula \eqref{eq3.4}.
\end{theorem}

If one considers from the start  the vector function $v(x)$ in \eqref{eq1.9} normalized by the condition $\|v\|=1$ (without loss in generality due to renormalization $\mu\rightarrow\mu\|v\|^2$), the number $\mu$ is unambiguously defined by the numbers $\mu|v_k|^2$ since in this case $\sum\limits_k|v_k|^2=1$. The function $v(x)$ is restored (ambiguously) by the formula
\begin{equation}
v(x)=\sum\limits_k|v_k|\theta_ku(\lambda_k(0),x)\label{eq3.6}
\end{equation}
where $\theta_k$ are arbitrary numbers such that $\theta_k\in\mathbb{T}$.

So from the spectrum $\sigma(D(\mu,v))$, one unambiguously defines the number $\mu$ and the class of functions $v(x)$ of the form \eqref{eq3.6} such that spectrum of the operator $D$
\eqref{eq1.9} coincides with $\sigma(D(\mu,v))$.
\vspace{5mm}

\section{Inverse problem. Inverse problem data description}\label{s4}

{\bf 4.1.} In Section \ref{s3}, method of ambiguous restoration of the function $v(x)$ from the spectral data is presented. Hereinafter a technique of unambiguous finding of $v(x)$ is given, though,
from the two spectra. The following statement belongs to B. V. Hvedelidze \cite{7}.

\begin{theorem}\label{t4.1}
{\bf \cite{7}}. If $\varphi(x)\in L^p(\mathbb{R})$ ($p>1$), then the Cauchy type integral
\begin{equation}
\Phi(z)=\frac1{\pi i}\int\limits_\mathbb{R}\frac{\varphi(x)dx}{x-z}\label{eq4.1}
\end{equation}
has non-tangential values a. e. on $\mathbb{R}$ from the half-planes $\mathbb{C}_\pm$ and the Sokhotsky formulas are true,
\begin{equation}
\Phi_\pm(x)=\Phi(x\pm i0)=\pm\varphi(x)+\frac1{\pi i}\not\hspace{-2.25mm}\int\limits_\mathbb{R} \frac{\varphi(t)}{t-x}dt,\label{eq4.2}
\end{equation}
besides, $\Phi_\pm(x)\in L^p(\mathbb{R})$.
\end{theorem}

The function $e^{i(\lambda\pi+\beta-\alpha)}R(\lambda)$ ($R(\lambda)$ is given by \eqref{eq1.28}) as $\lambda\in\mathbb{R}$ belongs to $L^2(\mathbb{R})$ due to the Paley -- Wiener theorem \cite{4}, therefore $\Re\{e^{i(\lambda\pi+\beta-\alpha)}R(\lambda)\}\in L^2(\mathbb{R})$ as $\lambda\in\mathbb{R}$. Limit values of the Cauchy type integral $\Phi(z)$ \eqref{eq4.1} from $\mathbb{C}_+$ on $\mathbb{R}$  when $\varphi(x)=\Re\{e^{i(\pi x+\beta-\alpha)}R(x)\}$ due to \eqref{eq4.2} are
$$\Phi_+(x)=\Re\{e^{i(\pi x+\beta-\alpha)}R(x)\}+\frac1{\pi i}\hspace{1.7mm}/\hspace{-4.25mm}\int\limits_\mathbb{R}\frac{\Re\{e^{i(\pi t+\beta-\alpha)}R(t)\}}{t-x}dt,$$
besides, integral in this sum equals $\Im\{e^{i(\pi x+\beta-\alpha)}R(x)\}$. Thus,
$$\Phi_+(x)=e^{i(\pi x+\beta-\alpha)}R(x)$$
and, after an analytical extension into $\mathbb{C}_+$,
\begin{equation}
e^{i(\lambda\pi+\beta-\alpha)}R(\lambda)=\Phi(\lambda)\quad(\lambda\in\mathbb{C}_+).\label{eq4.3}
\end{equation}
And since, due to \eqref{eq1.27},
$$\Re\{e^{i(\pi x+\beta-\alpha)}R(x)\}=-\frac2\mu(\Delta(\mu,x)-\Delta(0,x))\quad(x\in\mathbb{R}),$$
then \eqref{eq4.3} implies
$$e^{i(\lambda\pi+\beta-\alpha)}R(\lambda)=\frac{2i}{\mu\pi}\int\limits_\mathbb{R}\frac{\Delta(\mu,x)-\Delta(0,x)}{x-\lambda}dx\quad(\lambda\in\mathbb{C}_+)$$
or
\begin{equation}
R(\lambda)=\frac{2i}{\mu\pi}e^{-i(\lambda\pi+\beta-\alpha)}\int\limits_\mathbb{R}\frac{\Delta(\mu,x)-\Delta(0,x)}{x-\lambda}dx\quad(\lambda\in\mathbb{C}_+).\label{eq4.4}
\end{equation}
So, the function $R(\lambda)$ is unambiguously defined by the characteristic function $\Delta(\mu,\lambda)$ via formula \eqref{eq4.4}.

For the function $R(\lambda)$ \eqref{eq1.28}, due to \eqref{eq1.26}, the identiity
$$R(\lambda)+R^*(\lambda)=[e^{i\alpha}\widetilde{v}_+(\lambda)+e^{-i\alpha}\widetilde{v}_-(-\lambda)][e^{i\alpha}\widetilde{v}_+(\lambda)+e^{-i\alpha}\widetilde{v}_-(-\lambda)]^*$$
is true and, as $\lambda\in\mathbb{R}$,
\begin{equation}
2\Re R(\lambda)=|e^{i\alpha}\widetilde{v}_+(\lambda)+e^{-i\alpha}\widetilde{v}_-(-\lambda)|^2=|V(\lambda)|^2\quad(\lambda\in\mathbb{R})\label{eq4.5}
\end{equation}
where
\begin{equation}
V(\lambda)\stackrel{\rm def}{=}e^{i\alpha}\widetilde{v}_+(\lambda)+e^{-i\alpha}\widetilde{v}_-(-\lambda).\label{eq4.6}
\end{equation}
Formula \eqref{eq1.19} implies that
\begin{equation}\label{eq4.7}
\begin{array}{ccc}
{\displaystyle V(\lambda)=e^{i\alpha}\int\limits_0^\pi e^{-i\lambda x}(v_1(x)+iv_2(x))dx+e^{-i\alpha}\int\limits_0^\pi e^{i\lambda x}(v_1(x)-}\\
{\displaystyle -iv_2(x))dx=2\int\limits_0^\pi[\cos(\lambda x-\alpha)v_1(x)+\sin(\lambda x-\alpha)v_2(x)]dx}
\end{array}
\end{equation}

\begin{lemma}\label{l4.1}
The functions $v_1(x)$ and $v_2(x)$ are expressed via the function $V(\lambda)$ \eqref{eq4.6} by the formulas
\begin{equation}
v_1(x)=\frac1{2\pi}\int\limits_\mathbb{R}V(\lambda)\cos(\lambda x-\alpha)dx;\quad v_2(x)=\frac1{2\pi}\int\limits_\mathbb{R}V(\lambda)\sin(\lambda x-\alpha)dx.\label{eq4.8}
\end{equation}
\end{lemma}

P r o o f. Formula \eqref{eq4.7} implies that
$$V(-\lambda)=2\int\limits_0^\pi[\cos(\lambda x+\alpha)v_2(x)]dx,$$
therefore for
$$V_\pm(\lambda)=\frac12(V(\lambda)\pm V(-\lambda))$$
the representations
$$V_+(\lambda)=2\int\limits_0^\pi\cos\lambda x[v_1(x)\cos\alpha-v_2(x)\sin\alpha]dx;$$
$$V_-(\lambda)=2\int\limits_0^\pi\sin\lambda x[v_1(x)\sin\alpha+v_2(x)\cos\alpha]dx$$
are true. Applying the inverse sine and cosine Fourier transforms \cite{4}, one finds
$$v_1(x)\cos\alpha-v_2(x)\sin\alpha=\frac1\pi\int\limits_0^\infty\cos\lambda xV_+(\lambda)d\lambda;$$
$$v_1(x)\sin\alpha+v_2(x)\cos\alpha=\frac1\pi\int\limits_0^\infty\sin\lambda xV_-(\lambda)d\lambda,$$
therefore
$$v_1(x)=\frac1\pi\int\limits_0^\infty[\cos\lambda x\cos\alpha V_+(\lambda)+\sin\lambda x\sin\alpha V_-(\lambda)]d\lambda;$$
$$v_2(x)=\frac1\pi\int\limits_0^\infty[-\cos\lambda x\sin\alpha V_+(\lambda)+\sin\lambda x\cos\alpha V_-(\lambda)]d\lambda;$$
or
$$v_1(x)=\frac1{2\pi}\int\limits_0^\infty(\cos(\lambda x-\alpha)V(\lambda)+\cos(\lambda x+\alpha)V(-\lambda))d\lambda;$$
$$v_2(x)=\frac1{2\pi}\int\limits_0^\infty(\sin(\lambda x-\alpha)V(\lambda)-\sin(\lambda x+\alpha)V(-\lambda))d\lambda.$$

\begin{remark}\label{r4.1}
In the case when the vector function $v(x)=\col[c,c]$ ($c$ is a real number) is constant, function $V(\lambda)$ \eqref{eq4.7} equals
\begin{equation}
V(\lambda)=C(\lambda)=\frac{c4\sqrt2}\lambda\sin\frac{\lambda\pi}2\cdot\sin\left(\lambda\pi-\alpha+\frac\pi4\right).\label{eq4.9}
\end{equation}
\end{remark}

\begin{theorem}\label{t4.2}
By the two spectra $\sigma(D(\mu,v))$ and $\sigma(D(\mu,v+c))$ of the operators $D(\mu,v)$ and $D(\mu,v+c)$ \eqref{eq1.9} ($c=\col[c,c]$, $c\in\mathbb{R}$ and $c\not=0$) the number $\mu$ and real vector function $v(x)$ are unambiguously found.
\end{theorem}

P r o o f. Characteristic function $\Delta(\mu,\lambda)$ is unambiguously defined by its zeros $\lambda_k(\mu)$ that form the set $\sigma(D(\mu,v))$ (see Sec. \ref{s3}). From formula \eqref{eq4.4} one obtains the function $\mu R(\lambda)$ and from it, due to \eqref{eq4.5}, calculates the function $\mu V^2(\lambda)$. Similarly, by the second spectrum $\sigma(D(\mu,v+c))$ one finds
$$\mu(V(\lambda)+C(\lambda))^2=\mu(V^2(\lambda)+2V(\lambda)C(\lambda)+C^2(\lambda))$$
where $C(\lambda)$ is given by \eqref{eq4.9}. Knowing $\mu V^2(\lambda)$ and $\mu(V(\lambda)+C(\lambda))^2$, one defines $\mu V(\lambda)$. Using Lemma \ref{l4.1}, one calculates $\mu v_1(x)$, $\mu v_2(x)$. From normalization $\|v(x)\|=1$ one finds the number $\mu^2$, and sign of the number $\mu$, according to Remark \ref{r2.1}, is found depending on the nature of intermittency of numbers $\lambda_k(\mu)\in\sigma(D(\mu,v))$ and $\sigma_1$ \eqref{eq2.6}. $\blacksquare$
\vspace{5mm}

{\bf 4.2.} Present a more simple technique of recovery of the number $\mu$ and vector function $v(x)$ from the two spectral data. Fourier coefficients $v_k=\langle v(x), u(\lambda_k(0),x)$ are
$$v_k=\frac1{\sqrt{\pi}}\int\limits_0^\pi [v_1(x)\cos(\lambda_k(0)x-\alpha)+v_2(x)\sin(\lambda_k(0)x-\alpha)]dx$$
$$=\frac1{\sqrt\pi}\int\limits_0^\pi\{v_1(x)(\cos kx\cos(mx-\alpha)-\sin kx\sin(kx-\alpha))$$
$$+v_2(x)(\sin kx\cos(mx-\alpha)+\cos kx\sin(mx-\alpha))\}dx$$
where
$$\lambda_k(0)x-\alpha=kx+\frac{\alpha-\beta}\pi x-\alpha=kx+mx-\alpha\quad\left(m=\frac{\alpha-\beta}\pi\right).$$
Therefore
\begin{equation}
v_k=\frac1{\sqrt\pi}\int\limits_0^\pi\cos kxw_1(x)dx+\frac1{\sqrt\pi}\int\limits_0^\pi\sin kxw_2(x)dx\label{eq4.10}
\end{equation}
where
\begin{equation}
\begin{array}{ccc}
w_1(x)=v_1(x)\cos(mx-\alpha)+v_2(x)\sin(mx-\alpha);\\
w_2(x)=-v_1(x)\sin(mx-\alpha)+v_2(x)\cos(mx-\alpha).
\end{array}\label{eq4.11}
\end{equation}
Formula \eqref{eq4.10} implies
\begin{equation}
\frac{v_k+v_{-k}}2=\frac1{\sqrt\pi}\int\limits_0^\pi\cos kxw_1(x)dx;\frac{v_k-v_{-k}}2=\frac1{\sqrt\pi}\int\limits_0^\pi\sin kxw_2(x)dx.\label{eq4.12}
\end{equation}
Each function system $\left\{{\displaystyle \frac1{\sqrt\pi},\sqrt{\frac2\pi}\cos kx,k\in\mathbb{N}}\right\}$ and ${\displaystyle\left\{\sqrt\frac2\pi\sin kx,k\in\mathbb{N}\right\}}$ forms an orthonormal basis in $L^2[0,\pi]$ \cite{4}. Therefore the functions $w_1(x)$ and $w_2(x)$ are expanded into Fourier series by these bases and taking into account \eqref{eq4.12} one obtains
$$w_1(x)=\frac1{\sqrt\pi}v_0+\sum\limits_{k=1}^\infty(v_k+v_{-k})\frac1{\sqrt\pi}\cos kx;$$
$$w_2(x)=\sum\limits_{k=1}^\infty(v_k+v_{-k})\frac1{\sqrt\pi}\sin kx,$$
or, according to \eqref{eq4.11},
$$v_1(x)\cos(mx-\alpha)+v_2(x)\sin(mx-\alpha)=\frac1{\sqrt\pi}\left(v_0+\sum\limits_{k=1}^\infty(v_k+v_{-k})\cos kx\right);$$
$$-v_1(x)\sin(mx-\alpha)+v_2(x)\cos(mx-\alpha)=\frac1{\sqrt\pi}\sum\limits_{k=1}^\infty(v_k-v_{-k})\sin kx.$$
Hence it follows an analogue of Lemma \ref{l4.1}.

\begin{lemma}\label{l4.2}
By the Fourier coefficients $v_k=\langle v(x),u(\lambda_k(0),x)\rangle$ ($k\in\mathbb{Z}$, $u(\lambda_k(0),x)$ are given by \eqref{eq1.8}) the vector function $v(x)=\col(v_1(x),v_2(x))$ is calculated via the formulas
\begin{equation}
v_1(x)=\frac1{\sqrt\pi}\sum\limits_{-\infty}^\infty v_k\cos(\lambda_k(0)x-\alpha);\quad v_2(x)=\frac1{\sqrt\pi}\sum\limits_{-\infty}^\infty v_k\sin(\lambda_k(0)x-\alpha).\label{eq4.13}
\end{equation}
\end{lemma}

\begin{remark}\label{r4.2}
If $v(x)=c=\col[c,c]$ ($c\in\mathbb{R}$, $c\not=0$) is a constant vector function, then its Fourier coefficients $c_k=\langle c,u(\lambda_k(0),x)\rangle$ are
\begin{equation}
c_k=\frac{c2\sqrt2}{\sqrt\pi\lambda_k(0)}\sin\frac{\pi\lambda_k(0)}2-\sin\left(\frac{\pi\lambda_k(0)}2-\alpha+\frac\pi4\right)\quad(k\in\mathbb{Z}).\label{eq4.14}
\end{equation}
\end{remark}

Method of recovery of the number $\mu\in\mathbb{R}$ and {\bf real} vector function $v(x)$ by the two spectra consists in the following. According to Theorem \ref{t3.2}, the numbers $\mu v_k^2$ are calculated by formula \eqref{eq3.5}, and the numbers $\mu(v_k+c_k)^2=\mu(v_k^2+2v_kc_k+c_k^2)$ are calculated correspondingly by the spectrum $\sigma(D(\mu,v+c))$, $c_k$ are given by \eqref{eq4.14}. From these two sequences, $\mu v_k^2$ and $\mu(v_k+c_k)^2$, one finds $\mu v_k$ and in view of Lemma \ref{l4.2} via the formulas \eqref{eq4.13} defines $\mu v_1(x)$ and $\mu v_2(x)$, i. e., the vector function $\mu v(x)$. The number $\mu$ is defined from the normalization $\|v\|=1$ (see proof of Theorem \ref{t4.2}).

This method of calculation $\mu$ and $v(x)$ is preferable since it does not demand construction of the characteristic function $\Delta(\mu,\lambda)$ by the spectrum $\sigma(D(\mu,\lambda))$ and calculation of the singular integral \eqref{eq4.4}.
\vspace{5mm}

{\bf 4.3} Description of sequences of real numbers $\{\lambda_k(\mu)\}$ that form spectrum of the operator $D(\mu,v)$ \eqref{eq1.9} is more convenient to give in terms of description of the class of characteristic functions $\Delta(\mu,\lambda)$.

Let two discrete sets $A$ and $B$ be given, elements of which have finite multiplicity (number of repeats). The operation {\bf + union}
$$C=A\uplus B$$
is a union of sets $A$ and $B$ in which elements $a\in A\cap B$ are repeated in $C$ the number of times which is equal to the sum of multiplicities of $a\in A$ and $a\in B$.

$A$ and $B$ are said to be {\bf partially interlaced} if there exist such two splittings
$$A=A_0\cup A_1\quad(A_0\cap A_1=\emptyset);\quad B=B_0\cup B_1\quad(B_0\cap B_1=\emptyset)$$
that the subsets $A_1$ and $B_1$ are interlaced, besides, $B=A_0\uplus B_1$ and number of the elements of $A_0\cap B_1$ is no more than finite.

\begin{definition}\label{d4.1}
A function $\Delta(\lambda)$ is said to be of class $J(\lambda_k(0))$, where $\lambda_k(0)$ are given by \eqref{eq1.7}, if

${\rm (i)}$ $\Delta(\lambda)$ is an entire function of exponential type $\sigma\leq\pi$ bounded on the real axis, $\Delta(0)\not=0$, and $\Delta^*(\lambda)=\Delta(\lambda)$;

${\rm(ii)}$ zeros $\{a_k\}_{-\infty}^\infty$ of the function $\Delta(\lambda)$ are real and simple, except for finite number that has multiplicity $2$ and coincide with $\lambda_s(0)$;

${\rm(iii)}$ sets $\{a_k\}_{-\infty}^\infty$ and $\{\lambda_k(0)\}$ partially interlace;

${\rm (iv)}$ the function
\begin{equation}
F(\lambda)=\frac{\Delta(\lambda)}{\sin(\lambda\pi-\alpha+\beta)}\quad(\alpha,\beta\in\mathbb{R})\label{eq4.15}
\end{equation}
has the property
\begin{equation}
\lim\limits_{y\rightarrow\infty}F(iy)=1;\label{eq4.16}
\end{equation}

${\rm(v)}$ the limit
\begin{equation}
\lim\limits_{y\rightarrow\infty}y|F(iy)-1|<\infty\label{eq4.17}
\end{equation}
exists.
\end{definition}

The function $\Delta(\mu,\lambda)$ \eqref{eq1.27} is of class $J(\lambda_k(0))$ due to \eqref{eq2.15}.

\begin{theorem}
If $\Delta(\lambda)\in J(\lambda_k(0))$, then there exists an operator $D(\mu,v)$ \eqref{eq1.9}, characteristic function $\Delta(\mu,\lambda)$ of which coincides with $\Delta(\lambda)$, $\Delta(\mu,\lambda)=\Delta(\lambda)$.
\end{theorem}

P r o o f. (i) implies that $\Delta(\lambda)$ is of class $C$ and
$$\Delta(\lambda)=\Delta(0){\rm v. p.}\prod\left(1-\frac\lambda{a_k}\right),$$
and (ii) implies that the function $F(\lambda)$ \eqref{eq4.15} is a real meromorphic function, zeros and poles of which are interlacing. Due to the Krein theorem \cite{5}, this  means that $F(\lambda)$ maps $\mathbb{C}_+\rightarrow\mathbb{C}_+$ (or $\mathbb{C}_+\rightarrow\mathbb{C}_-$). According to Chebotarev's theorem \cite{5}, such meromorphic functions are given by
$$F(\lambda)=a\lambda+b-\frac{A_0}\lambda+\sum A_k\left(\frac1{\lambda_k(0)-\lambda}-\frac1{\lambda_k(0)}\right)$$
where $A_k\geq0$ ($\leq0$), $a\geq0$, $b\in\mathbb{R}$ and
$$\sum\limits_k\frac{A_k}{b_k^2}<\infty.$$
Then \eqref{eq4.16} implies that $a=0$ and $A_0=0$ (since $\sin(\beta-\alpha)\not=0$), besides,
$$b-\sum\frac{A_k}{\lambda_k(0)}=1.$$
So,
$$F(\lambda)=1+\sum\limits_k\frac{A_k}{\lambda_k(0)-\lambda},$$
besides, numbers $A_k$ all are of the same sign ($A_k\geq0$ or $A_k\leq0$). Condition \eqref{eq4.17} implies that
\begin{equation}
\sum\limits_kA_k=\mu<\infty,\label{eq4.18}
\end{equation}
then $|v_k|=A_k/\mu>0$ and by formula \eqref{eq3.6} one can construct the function $v(x)$. Define now the operator $D(\mu,v)$ \eqref{eq1.9} where $\mu$ is given by \eqref{eq4.18} and $v(x)$ \eqref{eq3.6}, then Theorem \ref{t2.2} implies that characteristic function of this operator is $\Delta(\mu,\lambda)=F(\lambda)\cdot\Delta(0,\lambda)$. Thus $\Delta(\mu,\lambda)=\Delta(\lambda)$. $\blacksquare$

\renewcommand{\refname}{ \rm 

\newpage

Dr. Zolotarev V. A.

B. Verkin Institute for Low Temperature Physics and Engineering
of the National Academy of Sciences of Ukraine\\
47 Nauky Ave., Kharkiv, 61103, Ukraine

Department of Higher Mathematics and Informatics, V. N. Karazin Kharkov National University \\
4 Svobody Sq, Kharkov, 61077,  Ukraine

E-mail:vazolotarev@gmail.com

\end{Large}
\end{document}